\documentclass[12pt,reqno]{amsart}

\DeclareMathOperator{\Br}{Br}%
\DeclareMathOperator{\chr}{char}%
\DeclareMathOperator{\Gal}{Gal}%
\DeclareMathOperator{\ind}{ind}%
\DeclareMathOperator{\inv}{inv}%

\usepackage{amssymb}
\usepackage{hyperref}

\begin{document}

\title[Cyclic algebras, Schur indices, norms, and Galois modules]{Cyclic algebras, Schur indices, norms, and Galois modules}

\author[J\'{a}n Min\'{a}\v{c}]{J\'an Min\'a\v{c}$^*$}
\address{Department of Mathematics, Middlesex College,
\ University of Western Ontario, London, Ontario \ N6A 5B7 \ CANADA}
\thanks{$^*$Research supported in part by Natural Sciences and
Engineering Research Council of Canada grant R0370A01.}
\email{minac@uwo.ca}

\author[Andrew Schultz]{Andrew Schultz}
\address{Department of Mathematics,  University of Illinois at Urbana-Champaign, 1409 W. Green Street, Urbana, IL \ 61801 \ USA}
\email{acs@math.uiuc.edu}

\author[John Swallow]{John Swallow$^\dag$}
\address{Department of Mathematics, Davidson College, Box 7046, Davidson, North Carolina \ 28035-7046 \ USA}
\thanks{$^\dag$Supported in part by National Science
Foundation grant DMS-0600122.}
\email{joswallow@davidson.edu}

\begin{abstract}
    Let $p$ be a prime and suppose that $K/F$ is a cyclic extension of degree $p^n$ with group $G$, with $\chr F\neq p$.  Let $J$ be the $\mathbb{F}_p G$-module $K^\times/K^{\times p}$ of $p$th-power classes. In our previous paper we established precise conditions for $J$ to contain an indecomposable direct summand of dimension not a power of $p$. At most one such summand exists, and its dimension must be $p^i+1$ for some $0\le i<n$. We show that for all primes $p$ and all $0\le i<n$, there exists a field extension $K/F$ with a summand of dimension $p^i+1$. The main theme is the investigation of an interesting link between this dimension and the Schur index of a certain cyclotomic cyclic algebra.
\end{abstract}

\dedicatory{To Paulo Ribenboim}

\date{July 18, 2009; revised January 3, 2010}

\maketitle

\newtheorem{theorem}{Theorem}
\newtheorem*{theorem*}{Theorem}
\newtheorem{corollary}{Corollary}
\newtheorem*{proposition*}{Proposition}

\theoremstyle{definition}
\newtheorem{remark}{Remark}
\newtheorem{proposition}{Proposition}
\newtheorem{example}{Example}

\newcommand{\Fp}{\mathbb{F}_p}
\newcommand{\Q}{\mathbb{Q}}
\newcommand{\N}{\mathbb{N}}
\newcommand{\R}{\mathbb{R}}
\newcommand{\Z}{\mathbb{Z}}

\parskip=10pt plus 2pt minus 2pt

Central simple algebras, representation theory, and Galois theory are among the favorite topics of Paulo Ribenboim. With the infectiousness of his enthusiasm, Paulo conveyed the beauty of these topics to his students, friends, colleagues, and to the readers of his work. Here we show an example of the natural interactions between these topics.

Let $p$ be a prime and $K/F$ a cyclic extension of fields of degree $p^n$ with Galois group $G$.  Assume further that $\chr F\neq p$.  Let $K^\times$ be the multiplicative group of nonzero elements of $K$ and $J=J(K/F):=K^\times/K^{\times p}$ be the group of $p$th-power classes of $K$.  We see that $J$ is naturally an $\Fp G$-module. In our previous paper \cite{MSS1} we established the decomposition of $J$ into indecomposables, as follows.

For $i\in \N$ let $\xi_{p^i}$ denote a primitive $p^i$th root of unity, and for $0\le i\le n$ let $K_i/F$ be the subextension of degree $p^i$, with $G_i=\Gal(K_i/F)$.  We adopt the convention that for all $i$, $\{0\}$ is a free $\Fp G_i$-module.

\begin{theorem}[{\cite[Theorems 1, 2, and 3]{MSS1}}]\label{thm1}\

    Suppose
    \begin{itemize}
        \item $F$ does not contain a primitive $p$th root of unity or
        \item $p=2$, $n=1$, and $-1\notin N_{K/F}(K^\times)$.
    \end{itemize}
    Then
    \begin{equation*}
        J \simeq \bigoplus_{i=0}^n Y_i
    \end{equation*}
    where each $Y_i$ is a free $\Fp G_i$-module.

    Otherwise, let
    \begin{equation*}
        m=m(K/F):=
        \begin{cases}
            -\infty, & \xi_p\in N_{K/F}(K^\times), \\ \min \{s\colon
            \xi_p\in N_{K/K_s}(K^\times)\}-1, & \xi_p\notin
            N_{K/F}(K^\times).
        \end{cases}
    \end{equation*}
    Then
    \begin{equation*}
        J \simeq X \oplus \bigoplus_{i=0}^n Y_i
    \end{equation*}
    where $Y_i$ is a free $\Fp G_i$-module and $X$ is an indecomposable $\Fp G$-module of $\Fp$-dimension $p^m+1$ if $m\ge 0$ and $1$ if $m=-\infty$.
\end{theorem}

\noindent It is not difficult to show directly that the decomposition is unique. (Alternatively, see the result of Azumaya \cite[p.~144]{AnFu}.)

This theorem is an extension of the main theorems in \cite{MS1}, and these results have already been applied to Galois embedding problems in \cite{MS2} and \cite{MSS2}.

From the well-known result of Albert \cite{A} concerning embedding a cyclic extension of degree $p^i$ to a cyclic extension of degree $p^{i+1}$, we see that $\xi_p\in N_{K/K_s}(K^\times)$ for all $s\in\{0,1,\dots,n\}$ if $m(K/F)=-\infty$ and $\xi_p\in N_{K/ K_s}(K^\times)$ for all $s\in\{m+1,\dots,n\}$ if $m(K/F) >-\infty$. (In fact it is an easy exercise to show that $N_{K/K_i}(K^\times) \subseteq N_{K/K_j}(K^\times)$ if $i\le j$ directly without referring to Albert's result.)

The submodules $Y_i$ are produced naturally using norms from different layers of the tower of field extensions.  However, the remaining submodule $X$ is more mysterious, and we consider a first problem concerning the classification of all $\Fp G$-modules occurring as $J(K/F)$:

\noindent\hangindent\parindent\hangafter=0 Given $n\ge 1$ and $d$ an element of the set
\begin{equation*}
    \{1, p^0+1, \dots, p^{n-1}+1\},
\end{equation*}
does there exist a cyclic extension $K/F$ with $\xi_p\in F^\times$ such that the exceptional summand $X$ has dimension $d$?

\noindent (We remark that in the case of $n=1$, $\xi_p\in F^\times$, the full realization problem---realizing all possible isomorphism classes for the $\Fp G$-module $J(K/F)$---has been solved \cite{MS2}.)

It turns out that we may answer this question in the affirmative using a construction of cyclic division algebras due to Brauer-Rowen.  For a particular extension $K/F$ with $\Gal(K/F) = \langle \sigma \rangle \simeq \Z/p^n\Z$, the cyclic algebra $$A(K/F) = (K/F,\sigma,\xi_p)$$ will be investigated for its connection to the module-theoretic structure of $K^\times/K^{\times p}$.    For instance, we will see in section 1 that the invariant $m(K/F)$ above can be recast as
  \begin{equation}\label{eq:recast}
        m(K/F)=
        \begin{cases}
            -\infty, & \mbox{ if }A(K/F) \mbox{ splits}, \\ \min \{s:
                       K_s \mbox{ splits }A(K/F)\}-1, & \mbox{ otherwise}.
        \end{cases}
    \end{equation}

In section~\ref{se:algebras} we begin by recalling basic facts about cyclic algebras, including the algebra $A(K/F)$.  We go on to describe how the Schur index of $A(K/F)$ is bounded by $m(K/F)$.  In some important cases these two invariants coincide; however, we show in section \ref{se:invariants} that we can have a strict inequality.  Thus it appears that $m(K/F)$ is a genuinely new invariant of $K/F$.

Further investigation of $m(K/F)$ and its relation with the associated cyclotomic algebra $A(K/F)$ is of considerable interest. The case when this algebra is a division algebra is considered in the last section. We also highlight in section~\ref{se:intpart} the interesting relations between the integer part function and the basic rule for cyclic algebras when passing to a cyclic extension.  In section~\ref{se:local} we show that $m(K/F)\in \{-\infty,0\}$ in the case of local fields, and we provide examples for both cases. Finally, in section~\ref{se:bypass} we show how to use just a simple observation about norms to bypass cyclic algebras and get a direct proof of Theorem~\ref{thm2}. However, we stress the intrinsic interest of the relationship between $m(K/F)$ and the Schur index, and we hope that this paper will promote its further investigation.

\section{Cyclic Cyclotomic Algebras, $m(K/F)$, Schur's
Index, and Norms}\label{se:algebras}

First we reformulate $m(K/F)$ in terms of cyclic algebras and then we use the construction of Brauer-Rowen of suitable cyclic algebras.  We will then prove the following theorem:

\begin{theorem}\label{thm2}
    Let $n\in \N$ and $t\in \{-\infty, 0, 1, \dots, n-1\}$.
    Then there exists a cyclic extension $K/F$ of degree $p^n$
    with $\xi_p\in F^\times$ and $m(K/F)=t$.
\end{theorem}


Before turning to the proof of the theorem, we recall some basic facts about cyclic algebras. If $L/E$ is a cyclic extension of degree $r>1$, with Galois group $G=\Gal(L/E)=\langle \tau \rangle$, and $b\in E^\times$, then
\begin{equation*}
    B=(L/E,\tau,b)
\end{equation*}
is a central simple algebra such that
\begin{equation*}
    B=\bigoplus_{0\leq j<r}u^j L,
\end{equation*}
where $u^{-1}du=\tau(d)$ for all $d\in L$ and $u^r=b$. Thus $B$ is an $E$-algebra of dimension $r^2$ over $E$. We say that $\deg B:=r$. If $B\simeq M_s(D)$, the matrix algebra containing matrices of size $s\times s$ over some division algebra $D$, then we set $\ind B=\sqrt{\dim_E D}$. We call $\ind B$ the Schur index of $B$. We denote the order of $[B]$ in the Brauer group $\Br(E)$ by $\exp B$. Finally, we observe the following important connection:
\begin{equation}\label{eq:SplittingAndNorms}
    [B]=0\mbox{ in }\Br(E)\quad \Longleftrightarrow \quad b\in N_{L/E}(L^\times).
\end{equation}
In this case, we say that $B$ splits. For further details on cyclic
algebras we refer the reader to \cite[Chapter~15]{P} and
\cite[Chapter 7]{R}.

The particular cyclic algebra in which we will be most interested is
the cyclotomic cyclic algebra
\begin{equation*}
    A(K/F) := (K/F,\sigma,\xi_p),
\end{equation*}
where here $\Gal(K/F) = \langle \sigma \rangle \simeq \Z/p^n\Z$ and $\xi_p\in F$.  Notice that equivalence (\ref{eq:SplittingAndNorms}) serves as justification for equation (\ref{eq:recast}) from the introduction.

\begin{proof}[Proof of Theorem \ref{thm2}]
    We begin with a construction of Brauer-Rowen. (See \cite{Br} for the original construction and see \cite[\S~7.3]{R} and \cite[\S~6]{RT} for some nice variations of Brauer's construction.)

    First suppose $0\le t<n$.  Set $q=p^{n-t}$ and let $K=\Q(\xi_q)(\mu_1,\dots,\mu_{p^n})$, where $\xi_q$ is a primitive $q$th root of unity and the $\mu_i$ are indeterminates over $\Q$.  Observe that $K$ has an automorphism $\sigma$ of order $p^n$ fixing $\Q(\xi_q)$ and permuting the $\mu_i$ cyclically.

    Let $F=K^{\langle \sigma \rangle}$ be the subfield of $K$ fixed by $\langle \sigma \rangle$ and, for each $1\le i\le n$, let $K_i = K^{\langle \sigma^{p^i}\rangle}$.  Then $K/F$ is a cyclic extension of degree $p^n$ satisfying $\Q(\xi_p)\subseteq F$, and $G=\langle \sigma \rangle = \Gal(K/F)$.  Denote by $\bar\sigma$ the restriction of $\sigma$ to the subfield $K_{t+1}$.

    Now $A=A(K/F)$ is Brauer-equivalent to the cyclic algebra $R=(K_{t+1}/F, \bar\sigma, \xi_q)$ by \cite[Corollary~15.1b]{P}. (See also section~\ref{se:intpart} below.) On the other hand, the construction of Brauer-Rowen shows that $R$ is a division algebra of degree $p^{t+1}$ and exponent $p$  \cite[Theorem~7.3.8]{R}. Since $[A]=[R]\neq 0$, we have $\xi_p\not\in N_{K/F}(K^\times)$.

    For all $0\le i\le n$ we have
    \begin{equation*}
        (K/F, \sigma, \xi_p) \otimes_F K_{i} \simeq
        (K/K_i, \sigma^{p^i}, \xi_p)
    \end{equation*}
    by \cite[Lemma 6, p.~74]{D}. Therefore, since $K_{t+1}$ is a maximal subfield of $R$, $K_{t+1}$ splits $A$:
    \begin{equation*}
        [A\otimes_F K_{t+1}] = 0\in \Br(K_{t+1}).
    \end{equation*}
    Therefore $\xi_p\in N_{K/K_{t+1}}(K^\times)$.  Hence $m(K/F)\le t$.

    Suppose $m=m(K/F)<t$.  We have seen that $K_{m+1}/F$ splits $A$, but $m<t$ implies $p^{t+1} = \ind A \le [K_{m+1}:F] < p^{t+1}$, a contradiction (see
    \cite[13.4]{P}). Hence $m(K/F)=t$.

    Now suppose that $t=-\infty$.  Let $F$ be a number field containing $\xi_p$ and, if $p=2$, $\xi_4$.  Then the extension $F^c/F$ obtained by adjoining all $p$th-power roots of unity is the cyclotomic $\Z_p$-extension of $F$.  Let $K/F$ be the subextension of degree $p^n$ of $F^c/F$.  Then $G=\Gal(K/F)$ is cyclic and $K/F$ embeds in a cyclic extension of $F$ of degree $p^{n+1}$. Therefore $\xi_p\in N_{K/F}(K^\times)$, by a result of Albert \cite{A}, and hence $m(K/F)=-\infty$.
\end{proof}

\begin{example}[The Case $n=1$]
    We briefly describe how the case $n=1$ may be handled quite simply.

    For the case $m(K/F)=0$, we set $F=\Q(\xi_p)(X)$, where $X$ is a transcendental element over $\Q(\xi_p)$, and $K=F(\root{p}\of{X})$. Write $G=\Gal(K/F)$ as $\langle \sigma \rangle$ with $\sigma(\root{p}\of{X}) = \xi_p\root{p}\of{X}$. Then the cyclic algebra $A=A(K/F)$ is a symbol algebra $A=\left(\frac{X,\xi_p}{F,\xi_p}\right)$. (See, for instance, \cite[p.~284]{P}.) Furthermore,
    \begin{equation*}
        -[A]=\left[\left(\frac{\xi_p,X}{F,\xi_p}\right) \right]=[(E/F,\tau,X)] \in \Br(F),
    \end{equation*}
    where $E=F(\xi_{p^2})$ and $\tau(\xi_{p^2})=\xi_p\xi_{p^2}$. However, it is an easy exercise (solved in \cite[p.~380]{P}) that $[(E/F,\tau,X)]\neq 0$. Hence $A$ is not split, and $m(K/F)=0$ as required.

    The $m(K/F)=-\infty$ case follows as in the end of the proof of Theorem~\ref{thm2}.
\end{example}


\section{Local Fields}\label{se:local}

For extensions $K/F$ of local fields, we may deduce that $m(K/F)\in \{-\infty, 0\}$, confirming \cite{B}.  In fact, since we need only the main theorem of local class field theory (see \cite[Theorem~7.1]{Iw}), our conclusion that $m(K/F)\in \{-\infty,0\}$ is valid for all fields containing $\xi_p$ which satisfy this theorem.  A family of such fields was investigated by Neukirch and Perlis: fields $k$ such that the pair $(G_k,\bar k^\times)$ is class formation in the sense of Artin-Tate \cite[Chapter 14]{AT}.  In \cite[Theorem 2, p. 532]{NP} they prove that all such fields possess local class field theory.  In particular, they point out that the field $\Q_p^a$ of all algebraic numbers in $\Q_p$, and similarly the field $\R^a$ of all real algebraic numbers, are examples of such fields.

First, we approach the conclusion with a continuation of our consideration of cyclic algebras, again writing $A$ for the algebra $A(K/F)$ defined in section 1.  If $[A]= 0\in \Br(F)$, then $m(K/F)=-\infty$.  Otherwise, since $\ind A = \exp A$ for local fields (see \cite[Corollary 17.10b]{P}), the local invariant $\inv_F A$ of $A$ is $s/p$ with $s \in \N, p \nmid s$. Because
\begin{equation*}
    \inv_E(A\otimes_F E) = [E:F]\inv_F A
\end{equation*}
(see \cite[Proposition~17.10]{P}), we obtain that $\inv_{K_1}(A\otimes_F K_1) = 0$.  Hence $[A\otimes_F K_1]=0\in \Br(K_1)$ and $m(K/F)=0$, as desired.  (Indeed, this same argument shows that any exponent $p$ central simple algebra over a local field is split by every extension of degree $p$.)

Alternatively, we prove that $m(K/F)\in \{-\infty,0\}$ directly, without using cyclic algebras.  Suppose that $K/F$ is a cyclic extension of local fields of degree $p^n$. Assume that $\xi_p\in F$ and $\xi_p\notin N_{K/F}(K^\times)$. We want to show that $\xi_p\in N_{K/K_1} (K^\times)$.

Consider the natural map
\begin{equation*}
    F^\times/N_{K/F}(K^\times)\longrightarrow K^\times_1/
    N_{K/K_1}(K^\times),
\end{equation*}
induced by the inclusion $F^\times\hookrightarrow K^\times_1$.
(Recall that under the inclusion map, we have $N_{K/F}(K^\times)\subseteq N_{K/K_1}(K^\times)$.)  From \cite[Theorem~7.1]{Iw} we see that
\begin{equation*}
    F^\times/N_{K/F}(K^\times)\simeq \Gal(K/F)=C_{p^n}
\end{equation*}
and
\begin{equation*}
    K^\times_1/N_{K/K_1}(K^\times)=\Gal(K/K_1)\simeq C_{p^{n-1}}.
\end{equation*}

The natural map is therefore a group homomorphism with nontrivial kernel. Because $\xi_p N_{K/F}(K^\times)\in \frac{F^\times} {N_{K/F}(K^\times)}$ generates the smallest nontrivial subgroup of $\frac{F^\times}{N_{K/F}(K^\times)}$ we see that $\xi_p N_{K/F}(K^\times)$ belongs to the kernel of this map and $\xi_p\in N_{K/F}(K^\times/K_1^\times)$, as required.

Now suppose $F$ is a local field containing $\xi_{p^{n+1}}$ and suppose $a\in F^\times$ such that $a\notin F^{\times p}$
if $p>2$ or $a\notin -4F^{\times 4}$ if $p=2$. Let
$K=F(\root{p^n}\of{a})$.  Then we see that $K/F$ is a cyclic
extension of degree $p^n$. (See \cite[Theorems~8.1 and
9.1]{La}.) Furthermore, because $N_{K/F}(\xi_{p^{n+1}})=\xi_p$, we conclude that $m(K/F)=-\infty$.

\begin{example}[A Case of $m(K/F)=0$]
    Following \cite[\S~3]{HS} we provide a simple example
    of a local field extension with $m(K/F)=0$.

    Dirichlet's theorem shows that there exists a prime $q$ such that $q=1+p^n + kp^{n+1} = 1+p^n(1+kp)$ for some $k\in \N\cup \{0\}$.
For such $q$ we know that $\xi_{p^n} \in \Q_q$ and $\xi_{p^{n+1}} \notin \Q_q$ (see \cite[\S~3.3]{Ko}). On the other hand by \cite[\S~3.3]{Ko} we know that $\Q_q(\xi_q)$ is a totally ramified cyclic extension of $\Q_q$ of degree $q-1$. Let $K/\Q_q$ be the cyclic extension $K\subseteq \Q_q(\xi_q)$ of degree $p^n$. Then by \cite[Proposition~4.1]{HS} we see that $K$ cannot be embedded into a cyclic extension $L/\Q_q$ of degree $p^{n+1}$. Hence $m(K/\Q_q)=0$.
\end{example}

\section{Cocycles, Carrying, and Cyclic Algebras}\label{se:intpart}

In section 1 in the proof of Theorem \ref{thm2}, we used \cite[Corollary~15.1b]{P}:

\textit{Let $F\subseteq K\subseteq E$ be a tower of cyclic extensions such that $G=\Gal(E/F) = \langle \sigma \rangle$, $\bar{\sigma}=$ restriction of $\sigma$ on $K$, and $[E:F]=a=bq$ where $b=[K:F]$ and $q= \frac{a}{b}$. Then}
\begin{equation*}
    \left[\left(K/F,\bar{\sigma},\alpha\right)\right]=\left
    [\left(E/F,\sigma,\alpha^q\right)\right]\in\Br(F)
\end{equation*}
\textit{for each $\alpha\in F^\times$.}

In \cite{P} the proof of this corollary uses a case-by-case calculation of cocycles.  We feel that because of this, a beautiful and interesting connection between the integer parts of rational numbers and cohomology may be lost. We shall here highlight this connection.

For $r\in \Q$, let $[r]\in\Z$ denote the integer part of $r\in\Q$, that is, the value $[r]$ satisfying $[r]\le r< [r]+1$.  As in \cite[p.~278]{P} the essential part of the proof of this statement is to prove that the following two cocycles $\Psi_\alpha$ and $\Phi_{\alpha^q}$ determine the same cohomology class in $H^2(G,E^\times)$:
\begin{equation*}
    \Psi_\alpha(\sigma^i,\sigma^j)=
    \begin{cases}
        1, & k+l<b\\
        \alpha, & k+l \ge b,
    \end{cases}
\end{equation*}
where $i=b\left[\frac{i}{b}\right]+k$ and $j=b\left[\frac{j}{b}\right]+l$; and
\begin{equation*}
    \Phi_{\alpha^q}(\sigma^i,\sigma^j)=
    \begin{cases}
        1, &i+j < a \\
        \alpha^q, &i+j \ge a.
    \end{cases}
\end{equation*}
where $0 \le i,j < a$.

Let $M=\langle\alpha\rangle$ be the subgroup of $F^\times$ (and therefore also $E^\times$) generated by $\alpha$.  Then we show that our cocycles $\Psi_\alpha$ and $\Phi_{\alpha^q}$ already determine the same cohomology class in $H^2(G,M)$ and therefore by naturality also the same class in $H^2(G,E^\times)$.

Both classes in $H^2(G,M)$ determine the group extensions $\left[\Psi_\alpha\right]$ and $\left[\Phi_{\alpha^q} \right]$ of shapes
\begin{equation*}
    1 \longrightarrow M \longrightarrow H_\Psi \longrightarrow
    G \longrightarrow 1
\end{equation*}
and
\begin{equation*}
    1 \longrightarrow M \longrightarrow H_\Phi \longrightarrow
    G \longrightarrow 1.
\end{equation*}
Choose lifts $\hat \sigma \in H_\Psi$ and $\tilde \sigma \in H_\Phi$ for  $\sigma \in G$, and let $0 \leq i,j < a$ be given.  As above, write $i=b\left[\frac{i}{b}\right]+k$ and $j=b\left[\frac{j}{b}\right]+l$, and also write $i+j = a\left[\frac{i+j}{a}\right] + \overline{i+j}$.  Notice that if $k+l<b$ then we have both $\left[\frac{i}{b}\right]+\left[\frac{j}{b}\right] = \left[\frac{i+j}{b}\right]$ and also $\hat \sigma^i \hat \sigma^j = \hat \sigma^{\overline{i+j}}$, so $$\hat \sigma^i\alpha^{\left[\frac{i}{b}\right]}\hat \sigma^j \alpha^{\left[\frac{j}{b}\right]} = \hat \sigma^{\overline{i+j}}\alpha^{\left[\frac{i+j}{b}\right]}.$$  If, on the other hand, we have $k+l\geq b$, then $\left[\frac{i}{b}\right]+\left[\frac{j}{b}\right]+1 = \left[\frac{i+j}{b}\right]$; since $\hat \sigma^i \hat \sigma^j = \hat \sigma^{\overline{i+j}}\alpha$ by our cocycle condition, this gives $$\hat \sigma^i \alpha^{\left[\frac{i}{b}\right]}\hat \sigma^j \alpha^{\left[\frac{j}{b}\right]} = \hat \sigma^i \hat \sigma^j \alpha^{\left[\frac{i}{b}\right]+\left[\frac{j}{b}\right]} = \hat \sigma^{\overline{i+j}} \alpha^{\left[\frac{i}{b}\right]+\left[\frac{j}{b}\right]+1}= \hat \sigma^{\overline{i+j}}\alpha^{\left[\frac{i+j}{b}\right]}.$$

From this calculation we learn two things:  first, that the map $\gamma:H_\Phi \to H_\Psi$ defined by
\begin{align*}
    \gamma(\tilde{\sigma}^i) &= \hat \sigma^i\alpha^{\left[\frac{i}{b}
    \right]} \\ \gamma(\alpha) &= \alpha
 \end{align*}
is an isomorphism between $\left[\Phi_{\alpha^q}\right]$ and $\left[\Psi_\alpha\right],$ thus completing our proof; and, second, that the cocycle $\Psi_\alpha$ is essentially the difference $\left[\frac{i+j}{b}\right]-\left[\frac{i}{b} \right]-\left[\frac{j}{b}\right]$.  More specifically, this difference is the logarithm of $\Psi_\alpha$ with base $\alpha$. (That ``carrying" is a cocycle was already observed in \cite{I}.)  

The referee of this article introduced us to the following elegant and conceptual alternative proof of the main result of this section.  This proof also shows that this statement is a straightforward consequence of the bilinearity of the cup product.

Let $\Gamma$ be the Galois group of a separable closure $F_s$ of $F$ containing $E$.  Let $\chi \in H^1(\Gamma,\Q/\Z)$ be the character of $\Gamma$ given by $\chi(\gamma) = k/a$, where the restriction of $\gamma$ to $E$ is $\sigma^k$.  Let $\partial \chi$ be the image of $\chi$ under the connecting map $\partial$, $$\partial: H^1(\Gamma,\Q/\Z) \to H^2(\Gamma,\Z)$$ associated with $0 \to \Z \to \Q \to \Q/\Z \to 0.$ Viewing a given $\beta \in F^\times$ as an element in $H^0(\Gamma,F_s^\times)$, one can then check through explicit cocycle calculations that $$\partial \chi \cup \beta = [(E/F,\sigma,\beta)] \in \Br(F) \simeq H^2(\Gamma,F_s^\times).$$  In the case where $\beta = \alpha^q$, we therefore have $\partial \chi \cup \alpha^q = [(E/F,\sigma,\alpha^q)]$.  On the other hand, the cup product's bilinearity gives $\partial \chi \cup \alpha^q = \partial(q\chi) \cup \alpha$, and since $q\chi$ is the character of $\Gamma$ associated to $K/F$ (and its generator $\bar \sigma$) we therefore have $\partial (q\chi) \cup \alpha = [(K/F,\bar \sigma,\alpha)]$, as desired.

\section{Bypassing Cyclic Algebras with a Norm Calculation}\label{se:bypass}

Brauer-Rowen's construction of the cyclic algebras above suggests replacing $\Q(\xi_q)(\mu_1,\dots,\mu_{p^n})$ by any quotient field $L$ of a unique factorization domain $D$. Let $U(D)$ be the group of units of $D$. In this case we have
\begin{proposition}
    Let $G$ be a group of order $n$ acting on a unique factorization domain $D$ as a group of $n$ distinct automorphisms of $D$. Assume also that $U(D)$ is a $G$-invariant subgroup of $D$.  Then
    \begin{equation*}
        N_{L/L^G}(L)\cap U(D) = (U(D))^n.
    \end{equation*}
\end{proposition}

\begin{proof}
    If $u=w^n\in U(D)^n$ for some $w\in U(D)$, then $u=N_{L/L^G}(w)$ and therefore
    \begin{equation*}
        U(D)^n\subseteq N_{L/L^G}(L)\cap U(D).
    \end{equation*}
    Assume now that $u \in U(D)\cap N_{L/L^G}(L)$.  Then
    $u = N_{L/L^G}(w)$ where
    \begin{equation*}
        w = d\;\;\frac{p_1\dots p_t}{q_1\dots q_s}
    \end{equation*}
    and $d\in U(D)$ and $p_1\dots p_t,q_1\dots q_s$ are primes in $D$.

    Therefore
    \begin{equation*}
        \left(\prod_{g\in G} \prod^s_{i=1} g(q_i)\right)
        N_{L/L^G}(w)=d^n \prod_{g\in G} \prod^t_{i=1}g(p_i).
    \end{equation*}
    From the unique factorization property of $D$ we conclude
    that $t=s$ and for each $i \in \{1,\dots,t\}$ there exists
    $g\in G,j\in \{1,\dots,t\}$ and $d_i\in U(D)$ such that
    \begin{equation*}
        p_i=d_i g(q_j).
    \end{equation*}
    Hence $N_{L/L^G}(p_i)=d_i^n N_{L/L^G}(q_j)$ and we see that
    \begin{equation*}
        u = N_{L/L^G}(w) = d^n \prod_{i=1}^{t}d_i^n\in U(D)^n.
    \end{equation*}
\end{proof}

Using the proposition above, one can directly see, without the use of cyclotomic cyclic algebras, that the example employed in our original proof, namely
\begin{equation*}
    K/F = \Q(\xi_q)(\mu_1,\dots,\mu_{p^n}) / (\Q(\xi_q)(\mu_1,\dots, \mu_{p^n}))^{\langle \sigma \rangle},
\end{equation*}
with $q=p^{n-t}$, works.  In fact the following alternative theorem provides more examples. (We leave out the easy case when $t=-\infty$.)

\begin{theorem}\label{thm3}
    Let $\xi_{p^{n-t}}\in F$ but $\xi_{p^{n-t+1}}\notin F$, where $t\in \{0,1,\dots,n-1\}$. Let $C_{p^n}$ act on $L=F(\mu_1,\dots,\mu_{p^n})$ by extending a faithful linear representation $\psi:C_{p^n}\to GL_{p^n}(F)$ on an $F$-vector space $V\subseteq L$ generated by $\mu_1,\dots, \mu_{p^n}$. Let $E$ be a fixed field of $C_{p^n}$. Then $m(L/E)=t$.
\end{theorem}

\begin{proof}
    Set
    \begin{equation*}
        D=F[\mu_1,\mu_2,\dots,\mu_{p^n}].
    \end{equation*}
    Then $D$ is a unique factorization domain.
    Moreover, $L=F(\mu_1,\dots,\mu_{p^n})$ is the quotient field of $D$, and $C_{p^n}$ acts naturally on $L$ via an extension of the action from $V$ to $L$. Then by Galois theory we obtain a chain of fields
    \begin{equation*}
        E=L_0\subseteq L_1\subseteq\dots\subseteq L_n = L,
    \end{equation*}
    where $L_i$ is the fixed subfield of $L$ under the cyclic
    subgroup $C_i$ of order $p^{n-i}$. Then we deduce
    \begin{align*}
        \xi_p\in N_{L/L_i} &\Longleftrightarrow \xi_p \in D^{p^{n-i}}\\
        &\Longleftrightarrow \xi_{p^{n-i+1}}\in F \\
        &\Longleftrightarrow n-i+1 \le
        n-t \\
        &\Longleftrightarrow t+1\le i.
    \end{align*}
    Hence $m(L/E) = t$.
\end{proof}

\section{The Invariants $m$ and $\ind A$}\label{se:invariants}

The proof of Theorem~\ref{thm2} turns on the fact that for the particular extension $K/F$ we have $\ind A = p^{m+1}$.  It is interesting to ask whether this equality holds generally.  We show in this section that the answer is negative.

We begin
with the inequality
\begin{equation*}
    \ind A \le p^{m+1}
\end{equation*}
as follows.  Observe that by the definition of $m(K/F)$,
\begin{equation*}
    [A\otimes_F K_{m(K/F)+1}] = 0 \in \Br(K_{m(K/F)+1})
\end{equation*}
for $m(K/F)\neq -\infty$.  Hence the inequality holds in the case $m(K/F)\neq -\infty$.  The statement also holds for $m(K/F)=-\infty$, since $A$ splits if and only if $m(K/F)=-\infty$; in fact, in this case we obtain an equality.

We show that equality does not always hold by considering the following example in the number field case. Recall first that for number fields $\exp A=\ind A$. (See, for instance, \cite[Theorem~18.6]{P}.) Therefore $\ind A$ is either $1$ or $p$ since the exponent of $A$ divides $p$:
\begin{equation*}
    [\otimes^p A]=[(K/F,\sigma,1)] = 0 \in \Br(F).
\end{equation*}
Hence it is enough to produce a case when $m(K/F)>0$.

Let $p=2$, $c\in 4\Z\setminus \{0\}$, $a=1+c^2\notin \Z^2$, and $d\in\{1,-1\}$ such that $d(a+\sqrt{a})$ is not a sum of two squares in $\Q_2$. (For example, take $a=17$ and $d=-1$.) It is well-known that then
\begin{equation*}
    F=\Q<K_1=\Q(\sqrt{a})<K_2=\Q\left(\sqrt{d(a+
    \sqrt{a})}\right)
\end{equation*}
is a tower of fields with $K_2/F$ cyclic of order 4. (See
\cite[p.~33]{JLY}.)

Let $\hat K_i$, $i=1, 2$, denote the completion of $K_i$ with respect to any valuation $v$ on $K_i$ which extends the $2$-adic valuation on $\Q$.  Since $8\mid a-1$, we have $\hat K_1=\Q_2$ and then we may and do assume that $\hat K_1=\Q_2\subseteq\hat K_2$.

Since $d(a+ \sqrt{a})$ is not a sum of two squares in $\Q_2$, the quaternion algebra $(d(a+\sqrt{a}),-1)_{\Q_2}$ is nonsplit. Hence $-1\notin N_{\hat K_2/\Q_2}(\hat K_2)$ and therefore $-1\notin N_{K_2/K_1} (K_2^\times)$. (See \cite[p.~353]{P}.) We obtain then that $m(K/F)=1$.

\section{When $A$ is a Division Algebra}\label{se:when}

Observe that if $A$ is a division algebra, then $\ind A=p^n$ and the chain of inequalities
\begin{equation*}
    p^n=\ind A\leq p^{m+1}\leq p^n
\end{equation*}
forces the equality $\ind A=p^{m+1}$.  In this section we show how a natural construction gives additional field extensions $L_{n-k}/F_{n-k}$ with $\ind A = p^{m+1} = p^{k+1}$ for every $k<n$.

More precisely, we shall construct a tower of extensions
\begin{equation*}
    F=F_1\subseteq F_2 \subseteq \dots \subseteq F_n
\end{equation*}
and a corresponding tower of induced extensions $L_i:=KF_i$ with the following properties:  $[F_{i+1}:F_i]=p$, $i<n$; $[L_i:F_i]=p^n$; and $L_i/F_i$ is a cyclic extension with $\Gal(L_i/F_i)=\langle \sigma_i\rangle$.  Since the extensions $L_i$ and $F_{i+1}$ are linearly disjoint over $F_i$, the cyclotomic algebras $A_i=(L_i/F_i, \sigma_i, \xi_p)$ and $A_{i+1}=(L_{i+1}/F_{i+1}, \sigma_{i+1}, \xi_p)$ are related by $A_{i+1}=A_i\otimes_{F_{i}} F_{i+1}$.  (See \cite[Lemma~7, p.~74]{D}.)  Then, passing from $A_i$ to $A_{i+1}$, the index of $A_{i+1}$ will either remain that of $A_i$ or drop by a factor of $p$.  We shall choose $F_{i+1}$ in such a way that $\xi_p$ will be a norm from $L_{i+1}$ to a subfield whose codimension's $p$-factor is one larger than at the $i$th step. This will guarantee that $p\cdot \ind A_{i+1}=\ind A_i$, and we shall have full control of the behavior of $\ind A_i$ for each $i$.

\begin{proposition}\label{prop1}
    Suppose that $A$ is a division algebra.  Set
    $F_i=F(\xi_{p^i})$ and $L_i=K (\xi_{p^i})$ for each
    $i=1,2,\dots,n$. Further set $A_i=A \otimes_F F_i$.

    Then $L_i$ is linearly disjoint from $F_n$ over $F_i$ and
    \begin{equation*}
        \ind A_i = p^{m(L_i/F_i)+1} = p^{n-i+1}, \quad i=1,
        2,\dots,n.
    \end{equation*}
\end{proposition}

\begin{proof}
    We proceed by induction on $i$.  The base $i=1$ is simply the case $K_1/F_1=K/F$, which follows from the observation at the beginning of the section.  Hence we assume that $A$ is a division algebra and, for some $i\in \{1, 2, \dots, n-1\}$, we have $[L_i:F_i]=p^n$, $\ind A_i = p^{n-i+1}$, and $m(L_i/F_i) = n-i$.

    We claim that $\xi_{p^{i+1}}\notin L_i$.  Otherwise, since $F_i(\xi_{p^{i+1}})/F_i$ is an extension of degree $1$ or $p$, we deduce that $\xi_{p^{i+1}}\in T_i$, where $T_i$ is the subfield of $L_i/F_i$ with $[L_i:T_i]=p^i$. Without loss of generality we may assume that $\xi_{p^{i+1}}^{p^i}=\xi_p$. Then $\xi_p = N_{L_i/T_i} (\xi_{p^{i+1}})$, and we obtain $m(L_i/F_i)\leq n-i-1$, a contradiction.

    Hence $L_i$ and $F_{i+1}$ are linearly disjoint over $F_i$. Therefore $L_{i+1}/F_{i+1}$ is a Galois extension of degree $p^n$ and
    \begin{equation*}
        G=\Gal(L_i/F_i)\simeq\Gal(L_{i+1}/F_{i+1}).
    \end{equation*}

    Now choose $\sigma_{i+1}\in\Gal(L_{i+1}/F_{i+1})$ such that the restriction of $\sigma_{i+1}$ to $L_i$ is $\sigma_i$. (We assume that $\sigma_i$ is already defined by induction, where $\sigma_1=\sigma$.) Then we see that
    \begin{equation*}
        A_{i+1}=A_i\otimes_{F_i}F_{i+1}\simeq(L_{i+1}
        /F_{i+1},\sigma_{i+1},\xi_p).
    \end{equation*}
    We therefore obtain from \cite[Proposition~13.4(v)]{P} that
    \begin{equation*}
        \ind A_{i+1}\geq\frac{\ind A_i}{p}=p^{n-i}.
    \end{equation*}

    On the other hand, we show that $p^{m(L_{i+1}/F_{i+1})+1} \le p^{n-i}$, as follows.  Since $\xi_{p^{i+1}} \in F_{i+1}^\times$, we have
    \begin{equation*}
        \xi_{p} \in N_{L_{i+1}/T_{i+1}}(L_{i+1}^\times),
    \end{equation*}
    where $F_{i+1} \subseteq T_{i+1}\subseteq L_{i+1}$ and
    $[L_{i+1}:T_{i+1}]=p^{i}$.  Hence $m(L_{i+1}/F_{i+1})\le
    n-i-1$.

    Putting these last two equalities together with the equality of the previous section, we reach the following chain:
    \begin{equation*}
        \ind A_{i+1}\leq p^{m(L_{i+1}/F_{i+1})+1}
        \leq p^{n-i}\leq\ind A_{i+1}.
    \end{equation*}
    We obtain $m(L_{i+1}/F_{i+1}) = n-(i+1)$ and
    $p^{m(L_{i+1}/F_{i+1}) +1}=\ind A_{i+1}$, as desired.
\end{proof}

To include $m(K/F)=-\infty$ in the proposition, it is sufficient to continue the induction one step further.  Set $F_{n+1} = F(\xi_{p^{n+1}})$ and $L_{n+1}=K(\xi_{p^{n+1}})$.  Then again $L_{n+1}/F_{n+1}$ is a cyclic extension of degree $p^n$ and $A_{n+1}=A\otimes_F F_{n+1}$ splits. We conclude that $m(L_{n+1}/F_{n+1}) = -\infty$.

The referee of this article pointed out an alternative proof of Proposition \ref{prop1}.  Suppose $F_i \cap K = F$ is given (this is not difficult to see, especially when $p>2$), and write $A = \oplus_{0 \leq j <p^n}\ \mu^j K$, where $\mu^{p^n}=\xi_p$ and $\mu^{-1}k\mu = \sigma(k)$ for each $k \in K$.  Then we have $F_i := F(\xi_{p^i}) = F(\mu^{p^{n-i}}) \subset A$.  The centralizer of $F_i$ in $A$, written $C_A(F_i)$, is Brauer equivalent to $A_i$ (see \cite[Lemma 13.3]{P}).  Observe that $C_A(F_i)$ is a division algebra and hence $\ind C_A(F_i) = \root\of{\dim_{F_i} C_A(F_i)}$.  From the double centralizer theorem (\cite[Theorem 12.7(ii)]{P}) we see that $$\dim_F  C_A(F_i)  \cdot\dim_F F_i=\dim_F A.$$  We may write $$\dim_F  C_A(F_i) =  \dim_{F_i} C_A(F_i)\cdot \dim_F F_i, $$ and by using $\dim_F\ A = p^{2n}$ and assuming $\dim_F F_i = p^{i-1}$, we obtain $$\dim_{F_i} C_A(F_i) =  p^{2(n-i+1)}.$$  Hence we have $\ind A_i = p^{n-i+1}$, as required.

We have left to show that $\dim_F F_i = p^{i-1}$.  One can find that $x^{p^{i-1}}-\xi_p = 0$ is irreducible for odd $p$ in \cite[Chapter 6, Theorem 9.1]{La}, hence giving the necessary equality.  The case for $p=2$ is resolved in \cite[Lemma 15.4]{P}, where the polynomial $X^{2^n}+1$ is shown to be irreducible, so that $[F(\xi_{2^{n+1}}):F]=2^n$.  Since $[F(\xi_{2^{i+1}}):F(\xi_{2^i})]\leq 2$ for each $1 \leq i \leq n$, it therefore follows that $[F_i:F]=2^{i-1}$.  This gives rise to the following interesting

\begin{corollary}
Suppose that $n \geq 2$.  Then $\root\of{2} \not\in F$.
\end{corollary}

\begin{proof}
Suppose that $\root\of{2} \in F$. Then $$[F(\xi_8):F]=[F(\root\of{-1}):F]=2<2^{3-1},$$ a contradiction.
\end{proof}

\section{Acknowledgements}

We thank John Labute for his interest and for his valuable
remarks related to this paper.  We also thank the referee of this paper for providing us with elegant proofs in sections 3 and 6.  We would also like to thank the referee of a previous version of this paper who suggested the approach considered in section 4.

\end{document}